# ON A GENERALIZED FALSE DISCOVERY RATE

By Sanat K. Sarkar[1] and Wenge Guo

*Temple University and National Institute of Environmental Health Sciences*

The concept of $k$-FWER has received much attention lately as an appropriate error rate for multiple testing when one seeks to control at least $k$ false rejections, for some fixed $k \geq 1$. A less conservative notion, the $k$-FDR, has been introduced very recently by Sarkar [*Ann. Statist.* **34** (2006) 394–415], generalizing the false discovery rate of Benjamini and Hochberg [*J. Roy. Statist. Soc. Ser. B* **57** (1995) 289–300]. In this article, we bring newer insight to the $k$-FDR considering a mixture model involving independent $p$-values before motivating the developments of some new procedures that control it. We prove the $k$-FDR control of the proposed methods under a slightly weaker condition than in the mixture model. We provide numerical evidence of the proposed methods' superior power performance over some $k$-FWER and $k$-FDR methods. Finally, we apply our methods to a real data set.

**1. Introduction.** The classical idea of controlling at least one false discovery has been generalized recently to that of controlling at least $k$ false discoveries, for some fixed $k > 1$. The rationale behind it has been that often in practice one is willing to tolerate a few false rejections, so by controlling $k$ or more false rejections the ability of a procedure to detect more false null hypotheses can potentially be improved. The $k$-FWER, the probability of at least $k$ false rejections, is one such generalized error rate that has received considerable attention [5, 6, 12, 13, 15, 16, 19, 20, 27]. With $V_n$ and $R_n$ denoting, respectively, the total number of false rejections and the total number of rejections of null hypotheses in testing $n$ null hypotheses, it is defined as

$$\text{(1.1)} \qquad k\text{-FWER} = \Pr\{V_n \geq k\},$$

Received September 2007; revised April 2008.
[1]Supported by NSF Grant DMS-06-03868.

generalizing the traditional familywise error rate (FWER). Sarkar [19] has introduced the following alternative error rate generalizing the usual false discovery rate (FDR) of Benjamini and Hochberg [1]:

$$(1.2) \qquad k\text{-FDR} = E(k\text{-FDP}), \qquad \text{where } k\text{-FDP} = \frac{V_n I(V_n \geq k)}{R_n \vee 1},$$

with $I(A)$ denoting the indicator of the event $A$ and $R_n \vee 1 = \max(R_n, 1)$. It is the expected ratio of $k$ or more false rejections to all rejections of null hypotheses, and, as $k$-FDR $\leq k$-FWER, controlling it is a less conservative approach than controlling the $k$-FWER.

Given $p$-values corresponding to the null hypotheses, Sarkar [19] provided a stepup $k$-FDR procedure utilizing the $k$th order joint null distributions of the $p$-values. It was assumed that these $p$-values are either independent or positively dependent in a sense slightly stronger than assumed for $k = 1$ in proving the FDR control of the Benjamini–Hochberg (BH) procedure [3, 17]. Later, Sarkar and Guo [22] have given stepup as well as stepdown procedures based on the bivariate null distributions of the $p$-values, assuming the $p$-values are independent or positively dependent in the same sense as when $k = 1$.

Alternative $k$-FDR procedures with independent $p$-values are constructed in this article taking the approach of conservatively estimating the FDR for a fixed rejection region and using these estimates to produce FDR controlling procedures, as in [23, 24, 28]. For a single-step test with a nonrandom threshold, we derive a formula for the $k$-FDR of the test under the mixture model considered in [23] and many other subsequent papers. The formula offers a new insight into the notion of $k$-FDR in relation to that of the FDR. It provides a simple and intuitive upper bound to the $k$-FDR that can be thought of as a scaled version of the FDR, with the $(k-1)$-FWER in testing $n-1$ null hypotheses being the scale factor. Motivated by this, we consider conservative point estimates of the product of FDR and the probability of at least $k-1$ false rejections while testing $n-1$ null hypotheses, given a fixed rejection region for each null hypothesis. Then we develop through these estimates procedures (stepwise) that control the $k$-FDR at a given level.

One of the new $k$-FDR procedures developed is a generalized version of the BH FDR procedure, Procedure 1. Others are improved versions of Procedure 1 using a class of estimates of the number of true null hypotheses. The $k$-FDR control of these procedures is proved assuming that the $p$-values are independent with each having the $U(0, 1)$ distribution when the corresponding null hypothesis is true. This is a slightly weaker assumption than the i.i.d. mixture model.

The performances of our procedures are numerically compared with other relevant procedures. It is important to point out that while we are in the paradigm of controlling $k$ false rejections, $k$-FWER and other $k$-FDR procedures should be the relevant competitors. With that in mind, we numerically compare Procedure 1 with two $k$-FWER procedures in Sarkar [20] to see the extent of power improvement we have in a $k$-FDR procedure over a $k$-FWER procedure. This improvement is seen to be quite significant, especially for large number of hypotheses.

Considering an ideal situation where the number of true null hypotheses is given to us by an oracle, we determine the oracle procedure. It is a stepup procedure that mimics Procedure 1 with the number of true null hypotheses assumed known. We numerically compare the powers of different $k$-FDR procedures, those proposed here, Procedure 1 and its modification with a particular choice of the estimate of true null hypotheses, and the one in Sarkar [19], relative to the power of the oracle procedure.

Although this paper is motivated by the work of [23, 24], we have not fully pursued their line of research here. We keep our focus mainly on developing procedures controlling the $k$-FDR and not on estimating it. Furthermore, we have not taken the route of generalizing Storey's concept of positive FDR and the related $q$-value method. Finally, we obtain our results only in the finite sample setting.

The layout of the paper is as follows. The $k$-FDR formula under the mixture model is given in Section 2. Having briefly introduced in Section 3 a class of conservatively biased point estimates of the $k$-FDR based on this formula, we motivate our procedures controlling the $k$-FDR in Section 4. The findings of numerical studies are presented in Section 5. An application to a real data set is provided in Section 6. Some concluding remarks and additional numerical investigations are made in Section 7. Proofs of all the main results are given in the Appendix.

**2. The $k$-FDR under mixture model.** Given $n$ null hypotheses $H_1, \ldots, H_n$, consider testing if $H_i = 0$ (true) or $H_i = 1$ (false) simultaneously for $i = 1, \ldots, n$, based on their respective $p$-values $p_1, \ldots, p_n$. We first consider a single step multiple testing procedure rejecting each $H_i = 0$ if $p_i \leq t$ for some fixed, nonrandom $t \in (0,1)$ and derive a formula for the $k$-FDR of this procedure under the following model considered in [23].

*Mixture model*: Let $(p_i, H_i)$, $i = 1, \ldots, n$, be i.i.d. as $(p, H)$, where

$$\Pr(p \leq u | H) = (1-H)u + H F_1(u), \qquad u \in (0,1),$$
(2.1)
$$\Pr(H=0) = \pi_0 = 1 - \Pr(H=1),$$

for some cdf $F_1(u)$.

THEOREM 2.1. *Let*

(2.2) $$V_n(t) = \sum_{i=1}^{n} I(p_i \leq t, H_i = 0), \qquad R_n(t) = \sum_{i=1}^{n} I(p_i \leq t).$$

*Then, under the above mixture model, the k-FDR of the single-step test rejecting each $H_i = 0$ if $p_i \leq t$ is given by*

$$\text{(2.3)} \qquad k\text{-FDR}_n(t) = n\pi_0 t E\left\{\frac{I[V_{n-1}(t) \geq k-1]}{R_{n-1}(t)+1}\right\}.$$

The formula in Theorem 2.1 provides an insight into the $k$-FDR as a measure of generalized FDR as well as a direction toward developing procedures that control it. To see this, consider first $k = 1$ and notice that for the FDR the formula is given by

$$\text{(2.4)} \qquad \begin{aligned} \text{FDR}_n(t) &= n\pi_0 t E\left\{\frac{1}{R_{n-1}(t)+1}\right\} \\ &= E[V_n(t)]E\left\{\frac{1}{R_{n-1}(t)+1}\right\}. \end{aligned}$$

Of course, since $R_{n-1} \sim \text{Bin}[n-1, F(t)]$, with $F(t) = \pi_0 t + (1-\pi_0)F_1(t)$, we have

$$E\left\{\frac{1}{R_{n-1}(t)+1}\right\} = \frac{1-[1-F(t)]^n}{nF(t)} = \frac{\Pr\{R_n(t) \geq 1\}}{nF(t)};$$

that is, the formula (2.4) is same as the following alternative formula

$$\text{(2.5)} \qquad \begin{aligned} \text{FDR}_n(t) &= \frac{\pi_0 t}{F(t)}\Pr\{R_n(t) \geq 1\} \\ &= \frac{E[V_n(t)]}{E[R_n(t)]}\Pr\{R_n(t) \geq 1\}, \end{aligned}$$

given in Storey [23] and commonly used in many subsequent papers. Nevertheless, (2.4) offers a slightly different insight into the FDR than (2.5), and it is this insight that helps us in understanding what the $k$-FDR means as a generalization of the FDR. Writing the $k$-FDR as

$$\text{(2.6)} \qquad \begin{aligned} &k\text{-FDR}_n(t) \\ &= \Pr\{V_{n-1}(t) \geq k-1\}n\pi_0 t E\left\{\frac{1}{R_{n-1}(t)+1}\bigg|V_{n-1}(t) \geq k-1\right\}, \end{aligned}$$

we see that it is a combination of the $(k-1)$-FWER in testing $n-1$ null hypotheses and an FDR-type measure conditional on at least $k-1$ false rejections in testing $n-1$ null hypotheses.

For a fixed $t$, $I[V_{n-1}(t) \geq k-1]$ is a stochastically increasing function of $R_{n-1}(t)$, because $V_{n-1}(t)$ is so; whereas, $[R_{n-1}(t)+1]^{-1}$ is a decreasing function of $R_{n-1}(t)$. Using these, we get

$$\text{(2.7)} \qquad E\left\{\frac{1}{R_{n-1}(t)+1}\bigg|V_{n-1}(t) \geq k-1\right\} \leq E\left\{\frac{1}{R_{n-1}(t)+1}\right\}$$

(see Appendix A.4 for a proof). In other words, we have

$$(2.8) \qquad k\text{-FDR}_n(t) \leq \Pr\{V_{n-1}(t) \geq k-1\}\text{FDR}_n(t).$$

Storey [23] estimated $\text{FDR}_n(t)$, given a fixed rejection region $(0,t)$ for each null hypothesis, by using conservative point estimates of the quantity $\pi_0 t/F(t)$. Borrowing Storey's idea, we consider conservatively estimating the quantity

$$(2.9) \qquad \frac{\pi_0 t}{F(t)}\Pr\{V_{n-1}(t) \geq k-1\}$$

for estimating the $k$-FDR toward developing procedures that control it. Before we do that in the next section, it would be interesting to see what the quantity (2.9) means and how it is related to the original definition of the $k\text{-FDR}(t)$.

Since $V_{n-1} \sim \text{Bin}(n-1, \pi_0 t)$, the probability $\Pr\{V_{n-1}(t) \geq k-1\}$ is equal to $G(k-1, n-1, \pi_0 t)$, where

$$(2.10) \qquad G(k,n,u) = \sum_{j=k}^{n} \binom{n}{j} u^j (1-u)^{n-j}, \qquad 0 < u < 1.$$

Also, $E\{V_n(t)I(V_n(t) \geq k)\} = n\pi_0 t G(k-1, n-1, \pi_0 t)$ and $E\{R_n(t)I(R_n(t) \geq 1)\} = nF(t)$. Thus, we see that

$$(2.11) \qquad \frac{\pi_0 t}{F(t)}\Pr\{V_{n-1}(t) \geq k-1\} = \frac{E\{V_n(t)I(V_n(t) \geq k)\}}{E\{R_n(t)I(R_n(t) \geq 1)\}},$$

that is, the quantity (2.9) is the ratio of the expectations, conditional on at least one rejection, of the numerator and denominator terms in $k$-FDP, the expectation of which is the $k$-FDR. This is similar to what Storey [23] noted when $k=1$, that is, for the ratio $\pi_0 t/F(t)$. Storey also showed that $\pi_0 t/F(t)$ is the positive false discovery rate (pFDR) defined in [23] under the mixture model. Thus, the quantity (2.9) is also seen to be a combination of the $(k-1)$-FWER in testing $n-1$ null hypotheses and the pFDR. The right-hand side ratio in (2.11) when $k=1$ has been referred to as the marginal FDR (mFDR) in [26] where an optimal procedure controlling the mFDR is developed under the model (2.1), taking a compound-decision theoretic approach to multiple testing.

**3. Conservative point estimates of the $k\text{-FDR}(t)$.** Storey [23] proposes the following class of conservative point estimates of $\text{FDR}(t)$:

$$(3.1) \qquad \widehat{\text{FDR}}_\lambda(t) = \frac{n\hat{\pi}_0(\lambda)t}{R_n(t) \vee 1},$$

where

$$\hat{\pi}_0(\lambda) = \frac{n - R_n(\lambda)}{n(1 - \lambda)} \quad (3.2)$$

for any $\lambda \in [0, 1)$. Multiplying this with $G(k - 1, n - 1, t)$, a conservative version of $\Pr\{V_{n-1}(t) \geq k - 1\}$, we consider estimating the $k$-FDR$(t)$ as follows:

$$\widehat{k\text{-FDR}}_\lambda(t) = \frac{n\hat{\pi}_0(\lambda) t G(k - 1, n - 1, t)}{R_n(t) \vee 1}, \quad \lambda \in [0, 1). \quad (3.3)$$

THEOREM 3.1. *Let the p-values be independent and those corresponding to the true null hypotheses be i.i.d. $U(0, 1)$. Then, $E(\widehat{k\text{-FDR}}_\lambda(t)) \geq k\text{-FDR}(t)$, for every fixed $\lambda \in [0, 1)$.*

This result follows from [23, 24]. It shows that the point estimates given by (3.3) for the $k$-FDR are conservative.

REMARK 3.1. A more natural way of estimating the $k$-FDR$(t)$ would be to estimate $\Pr\{V_{n-1}(t) \geq k - 1\}$ using $G(k - 1, n - 1, \hat{\pi}_0(\lambda)t)$, instead of $G(k - 1, n - 1, t)$, and multiply this with (3.1). However, for such estimates, Theorem 3.1 holds under certain restrictions on $\lambda$ depending on $t$.

**4. Procedures controlling the $k$-FDR.** Using $\widehat{k\text{-FDR}}_\lambda(t)$, we will now derive a new class of $k$-FDR procedures. Let

$$t_\alpha(\widehat{k\text{-FDR}}_\lambda) = \sup\{0 \leq t \leq 1 : \widehat{k\text{-FDR}}_\lambda(t) \leq \alpha\}. \quad (4.1)$$

Then, reject $H_i$ if $p_i \leq t_\alpha(\widehat{k\text{-FDR}}_\lambda)$. Given $p_{1:n} \leq \cdots \leq p_{n:n}$, the sorted $p$-values, this procedure when $\lambda = 0$ (i.e., $\hat{\pi}_0 = 1$) is equivalent to the following procedure.

PROCEDURE 1. *Reject $H_{(1)}, \ldots, H_{(\hat{l})}$, where*

$$\hat{l} = \max\left\{1 \leq i \leq n : p_{i:n} G(k - 1, n - 1, p_{i:n}) \leq \frac{i\alpha}{n}\right\}, \quad (4.2)$$

*if the maximum exists, otherwise, reject none, where $H_{(i)}$ is the null hypothesis corresponding to $p_{i:n}$, $i = 1, \ldots, n$.*

THEOREM 4.1. *Procedure 1 controls the $k$-FDR at $\alpha$ if the p-values are independent and those corresponding to the true null hypotheses are i.i.d. $U(0, 1)$.*

Define

(4.3) $$\tilde{G}_{k,n}(t) = tG(k-1, n-1, t), \qquad t \in (0,1).$$

Let $\tilde{G}_{k,n}^{-1}$ be the inverse function of $\tilde{G}_{k,n}$. Then, Procedure 1 is a stepup procedure with the critical values $\alpha_i = \tilde{G}_{k,n}^{-1}(i\alpha/n)$, $i = 1, \ldots, n$, generalizing the BH procedure from an FDR to a $k$-FDR procedure. As $\tilde{G}_{k,n}(t) \leq t$ and $\tilde{G}_{k,n}(t)$ is increasing in $t$ (see, e.g., Result A.1 in Appendix A.2), $\tilde{G}_{k,n}^{-1}(t) \geq t$. In other words, Procedure 1 is uniformly more powerful than the BH procedure.

It is important to note that, as in a $k$-FWER procedure, the first $k-1$ critical values in Procedure 1 and the one to be developed later can be chosen arbitrarily. This is because the first $k-1$ critical values in any stepwise procedure have no role in defining the $k$-FDR of such a procedure as the $k$-FDR is zero until at least $k$ of the null hypotheses are rejected. Nevertheless, the best way to choose these critical values would be to keep them all constant at the $k$th critical value; see, for example, [19]. So, we consider the first $k-1$ critical values in our proposed $k$-FDR methods to be same as the $k$th one while comparing them with $k$-FWER and other $k$-FDR procedures.

Does Procedure 1 (with its first $k-1$ critical values same as the $k$th one) provide a more powerful method of controlling $k$ false rejections than a compatible $k$-FWER method? Lehmann and Romano [13] gave a stepdown $k$-FWER procedure generalizing Holm's original FWER procedure in [9]. Sarkar [20] showed that a stepup version of the procedure, which generalizes Hochberg's procedure in [8], also controls the $k$-FWER under independence or certain type of positive dependence. Its critical values are $\alpha_i = k\alpha/(n - i \vee k + k), i = 1, \ldots, n$. Clearly, Procedure 1 is more powerful than this procedure, as $i\alpha/n \geq k\alpha/(n+k-i)$, for all $i = k, \ldots, n$.

Theorem 4.1 establishes the $k$-FDR control of the single step procedure with the random threshold $t_\alpha(\widehat{k\text{-FDR}}_{\lambda=0})$. For $\lambda > 0$, we use the threshold, as in [24], based on the following modified version of $\widehat{k\text{-FDR}}_\lambda(t)$:

(4.4) $$\widehat{k\text{-FDR}}_\lambda^*(t) = \begin{cases} \dfrac{n\hat{\pi}_0^*(\lambda)tG(k-1, n-1, t)}{R_n(t) \vee 1}, & \text{if } t \leq \lambda, \\ 1, & \text{if } t > \lambda, \end{cases}$$

where

(4.5) $$\hat{\pi}_0^*(\lambda) = \frac{n - R_n(\lambda) + 1}{n(1-\lambda)},$$

a slight modification of Storey's [23] original estimate in (3.2). Define

(4.6) $$t_\alpha(\widehat{k\text{-FDR}}_\lambda^*) = \sup\{0 \leq t \leq 1 : \widehat{k\text{-FDR}}_\lambda^*(t) \leq \alpha\},$$

and reject $H_i$ if $p_i \leq t_\alpha(\widehat{k\text{-FDR}}_\lambda^*)$. Given $p_{1:n} \leq \cdots \leq p_{n:n}$, this is equivalent to finding $j = \max\{1 \leq i \leq n : p_{i:n} \leq \lambda\}$, for a fixed $\lambda \in (0,1)$, and then rejecting $H_{(1)}, \ldots, H_{(\hat{l})}$, where

$$(4.7) \qquad \hat{l} = \max\left\{1 \leq i \leq j : p_{i:n} \leq \min\left[\tilde{G}_{k,n}^{-1}\left(\frac{i\alpha(1-\lambda)}{n-j+1}\right), \lambda\right]\right\},$$

if the maximums at both stages exist, otherwise not rejecting any hypothesis. Nevertheless, we will consider slightly more conservative procedures of the following type.

PROCEDURE 2. Given a fixed $\lambda \in (0,1)$, find, at the first stage, $j = \max\{1 \leq i \leq n : p_{i:n} \leq \lambda\}$. At the second stage, reject $H_{(1)}, \ldots, H_{(\hat{l})}$, where

$$(4.8) \qquad \hat{l} = \max\left\{1 \leq i \leq j : p_{i:n} \leq \lambda \min\left[\tilde{G}_{k,n}^{-1}\left(\frac{i\alpha(1-\lambda)}{\lambda(n-j+1)}\right), 1\right]\right\}.$$

If the maximum does not exist at either stage, do not reject any hypothesis.

THEOREM 4.2. *Procedure 2, for every fixed $\lambda \in (0,1)$, controls the $k$-FDR at $\alpha$ if the $p$-values are independent and those corresponding to the true null hypotheses are i.i.d. $U(0,1)$.*

REMARK 4.1. Let $n_0$ be the number of true null hypotheses. When $n_0 < k$, the $k$-FDR is zero and hence trivially controlled. When $n_0 \geq k$, the $k$-FDR of a stepup procedure with critical values $\alpha_1 \leq \cdots \leq \alpha_n$ is bounded above by $n_0 \max_{k \leq r \leq n}\{\tilde{G}_{k,n_0}(\alpha_r)/r\}$ under the conditions assumed in Theorem 4.1, as seen from (A.2) and (A.8). With unknown $n_0$, we consider the maximum of this upper bound with respect to $n_0$, which is $n \max_{k \leq r \leq n}\{\tilde{G}_{k,n}(\alpha_r)/r\}$, and choose the $\alpha_r$ satisfying $n\tilde{G}_{k,n}(\alpha_r)/r = \alpha$ that makes it equal to $\alpha$. This is how we develop Procedure 1 generalizing the BH FDR procedure.

Procedure 2 is an adaptive $k$-FDR procedure generalizing that in [24]. It attempts to improve the $k$-FDR control of Procedure 1 by sharpening it using an estimate of $n_0$ obtained from the available $p$-values. More formally, given any $0 < \lambda < 1$, we choose the $\alpha_r$ subject to $\hat{n}_0(\lambda)\tilde{G}_{k,n}(\alpha_r)/r = \alpha$, considering the $\hat{n}_0(\lambda) = [n - R_n(\lambda) + 1]/(1-\lambda)$ used in [24], and then slightly modify it so that we can theoretically establish the $k$-FDR control of the modified adaptive procedure when $k > 1$. Of course, when $k = 1$ this modification does not make any difference, while it results in a slightly more conservative procedure when $k > 1$. We will explain later why it is more conservative when $k > 1$.

A more reasonable approach to constructing a class of adaptive procedures would be to find the $\alpha_r$ initially from $\hat{n}_0(\lambda)\tilde{G}_{k,\hat{n}_0(\lambda)}(\alpha_r)/r = \alpha$ before it is modified if necessary. We have taken a slightly more conservative approach than this, and the only reason we have done so is that we are able to theoretically prove the $k$-FDR control of the resulting adaptive procedure based on (4.8), but not of the original one based on (4.7) when $k > 1$. This proof is an extension of a proof given in [21] and alternative to those given in [2, 24] of the FDR control of the procedure in [24].

A careful study of our proof of the $k$-FDR control of Procedure 2 would reveal, at least theoretically, that (4.5) is a natural choice for an estimate of $n_0$ that can be used adaptively in Procedure 1 maintaining a control of the $k$-FDR. Alternative and more complicated methods of estimating $n_0$ are given in [10, 14]. However, unlike (4.5), it seems hard to prove theoretically that using these estimates adaptively in Procedure 1 will control the $k$-FDR. Of course, if these are used nonadaptively, that is, by obtaining them independently before incorporating into Procedure 1 to find $\alpha_r$ satisfying $\hat{n}_0 \tilde{G}_{k,n}(\alpha_r)/r = \alpha$, the $k$-FDR can be controlled as long as $E(1/\hat{n}_0) \le 1/n_0$.

If $n_0 \ge k$ were known, the least conservative stepup procedure controlling the $k$-FDR at level $\alpha$ would be the one in which $\alpha_r$ satisfies $n_0 \tilde{G}_{k,n_0}(\alpha_r)/r = \alpha$. This will be referred to as the oracle $k$-FDR procedure in this article.

Let us now explain why Procedure 2 is more conservative when $k > 1$ than the one based on (4.7). For any $0 < \lambda < 1$ and $t > 0$ [in particular, for $t = i\alpha(1-\lambda)/(n-j+1)$], note that

$$\tilde{G}_{k,n}(\lambda t) = \lambda t G(k-1, n-1, \lambda t) \le \lambda t G(k-1, n-1, t) = \lambda \tilde{G}_{k,n}(t).$$

Thus,

$$\lambda \tilde{G}_{k,n}^{-1}\left(\frac{t}{\lambda}\right) \le \tilde{G}_{k,n}^{-1}\left(\lambda \tilde{G}_{k,n}\left(\tilde{G}_{k,n}^{-1}\left(\frac{t}{\lambda}\right)\right)\right) = \tilde{G}_{k,n}^{-1}(t).$$

**5. Numerical studies.** Sarkar [20] introduced two stepup procedures that control the $k$-FWER under independence. One of these we call the generalized Hochberg procedure, and is based on the following critical values:

$$(5.1) \qquad \alpha_i^{(1)} = \frac{k\alpha}{n - i \vee k + k}, \qquad i = 1, \ldots, n.$$

As said before, this is actually the stepup analog of the $k$-FWER procedure derived in [13] as a generalization of the Holm procedure [9]. The other procedure we call Sarkar's $k$-FWER procedure, and is based on the following critical values:

$$(5.2) \qquad \alpha_i^{(2)} = \left(\alpha \prod_{j=1}^{k} \frac{j}{n - i \vee k + j}\right)^{1/k}, \qquad i = 1, \ldots, n.$$

In addition, Sarkar [19], while introducing the notion of the $k$-FDR, proposed a stepup $k$-FDR procedure with the following critical values:

$$\text{(5.3)} \qquad \alpha_i^{(3)} = \left( \frac{i \vee k}{n} \alpha \prod_{j=1}^{k-1} \frac{j}{n - i \vee k + j} \right)^{1/k}, \qquad i = 1, \ldots, n.$$

It controls the $k$-FDR at $\alpha$ under independence. We call this procedure Sarkar's $k$-FDR procedure. The oracle $k$-FDR procedure is a stepup procedure with the following critical values:

$$\text{(5.4)} \qquad \alpha_i^{(4)} = \tilde{G}_{k,n_0}^{-1}\left( \frac{i \vee k}{n_0} \alpha \right), \qquad i = 1, \ldots, n,$$

with $n_0 \geq k$.

Numerical studies were conducted, first to get an idea of how powerful the notion of $k$-FDR is compared to that of the $k$-FWER. For that, we considered Procedure 1 and compared it with the above $k$-FWER procedures, the generalized Hochberg and Sarkar's $k$-FWER procedures, in terms of their critical values and average powers. Second, we wanted to compare the average powers of different $k$-FDR procedures, Procedures 1 and 2 and Sarkar's $k$-FDR procedure, relative to the oracle $k$-FDR procedure. To recall the definition of average power, it is the expected proportion of false nulls that are rejected.

Figure 1 presents a comparison among Procedure 1, labeled New SU, and the generalized Hochberg and Sarkar's $k$-FWER procedures, labeled GH SU and Sarkar SU, respectively. We plot in this figure the three sequences of constants described in (4.2), (5.1) and (5.2) for $(n,k) = (500,8), (1000,10)$, $(2000, 15)$ and $(5000, 20)$ and $\alpha = 0.05$. The critical values of Procedure 1 are seen to be uniformly much larger than those of the generalized Hochberg procedure and, except when $n$ is quite large, they are also larger than those of Sarkar's $k$-FWER procedure.

Figure 2 presents a comparison among the above three procedures in terms of simulated average power. We considered in this case $(n, k) = (100, 3)$, $(200, 5), (500, 8)$ and $(1000, 10)$, and $\alpha = 0.05$. Each simulated average power was obtained by: (i) generating $n$ independent normal random variables $N(\mu_i, 1), i = 1, \ldots, n$, with $n_1$ of the $n$ $\mu_i$'s being equal to 2 and the rest 0; (ii) applying Procedure 1 and the generalized Hochberg and Sarkar's $k$-FWER procedures to the generated data to test $H_i : \mu_i = 0$ against $K_i : \mu_i \neq 0$ simultaneously for $i = 1, \ldots, n$ at $\alpha = 0.05$; and (iii) repeating steps (i) and (ii) 1000 times before observing the proportion of the $n_1$ false $H_i$'s that are correctly declared significant. As seen in this figure, Procedure 1 is uniformly much more powerful than the generalized Hochberg procedure and substantially more powerful than Sarkar's $k$-FWER procedure, with the power difference getting significantly higher with increasing number of false null hypotheses.

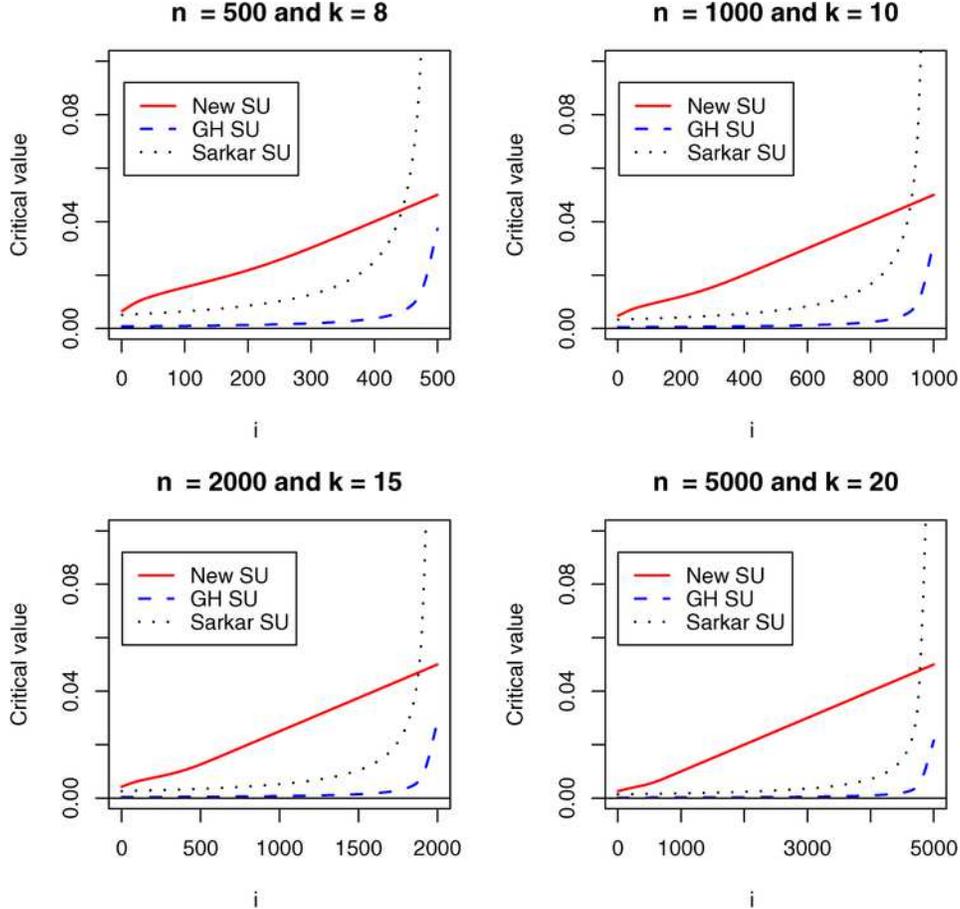

Fig. 1. *The critical constants of Procedure 1 (4.2), generalized Hochberg k-FWER procedure (5.1) and Sarkar's k-FWER procedure (5.2) for $\alpha = 0.05$.*

We brought in the above power study the other three $k$-FDR procedures, Procedure 2 (with $\lambda = 0.5$), Sarkar's $k$-FDR procedure and the oracle $k$-FDR procedure. Figure 3 presents this comparison, with Procedure 1 now labeled New SU I and Procedure 2, Sarkar's $k$-FDR procedure and the oracle procedure labeled, respectively, New SU II, Sarkar and Oracle. Benchmarking the three $k$-FDR procedures, Procedures 1 and 2 and Sarkar's procedure against the oracle, it is seen that Procedure 1 has the best power performance among these three when the number of false null hypotheses is small. But, with increasing number of false null hypotheses, Procedure 2 becomes substantially more powerful than either of the other two procedures.

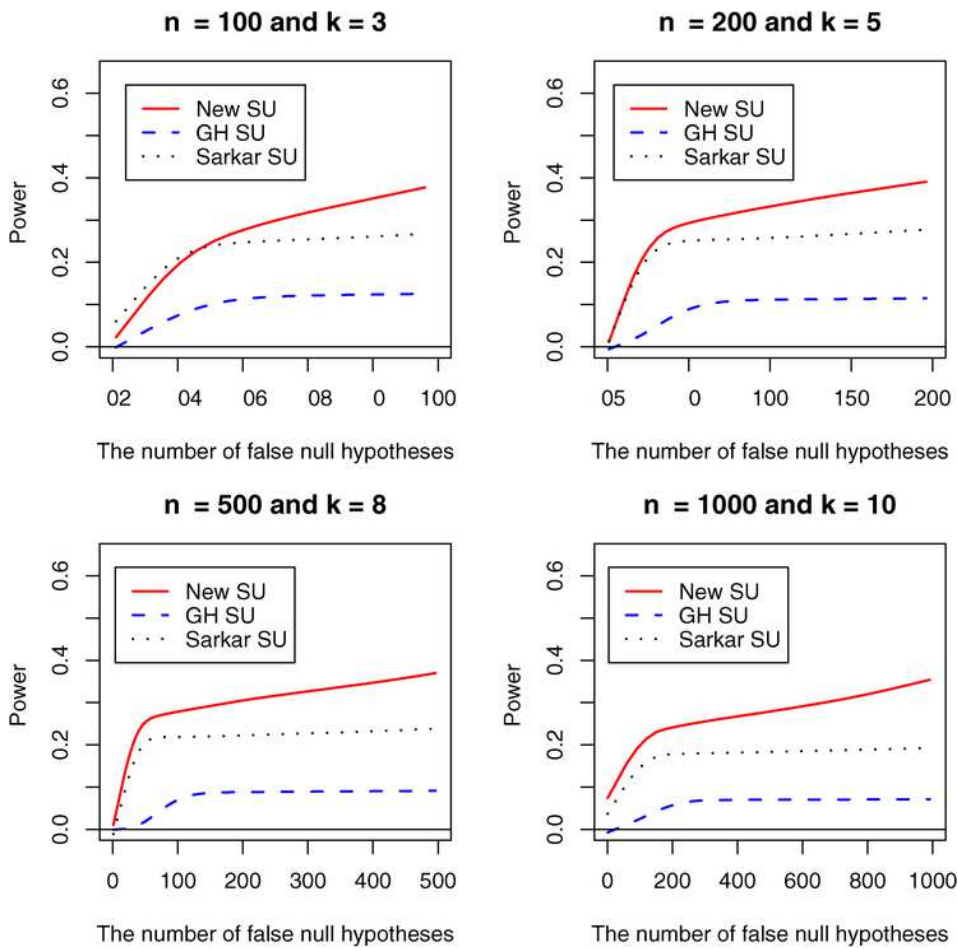

FIG. 2. *Comparison of average powers of k-FDR stepup procedures based on the sets of critical values given by (4.2), (5.1) and (5.2) for $\alpha = 0.05$.*

**6. An application to gene expression data.** Hereditary breast cancer is known to be associated with mutations in BRCA1 and BRCA2 proteins. Hedenfalk et al. [7] report a group of differentially expressed genes between tumors with BRCA1 mutations and tumors with BRCA2 mutations by analyzing one real microarray data set. The data set, which is publicly available from the web site http://research.nhgri.nih.gov/microarray/NEJM_Supplement/, consists of 22 breast cancer samples, among which 7 samples are BRCA1 mutants, 8 samples are BRCA2 mutants, and the remaining 7 samples are sporadic (not used in this illustration). Expression levels in terms of florescent intensity ratios of a tumor sample to a common reference sample, are measured for 3226 genes using cDNA microarrays. Before processing the data, there is a preprocessing step. If any gene has one ratio exceeding 20, then this gene is eliminated. Such preprocessing leaves $n = 3170$ genes.

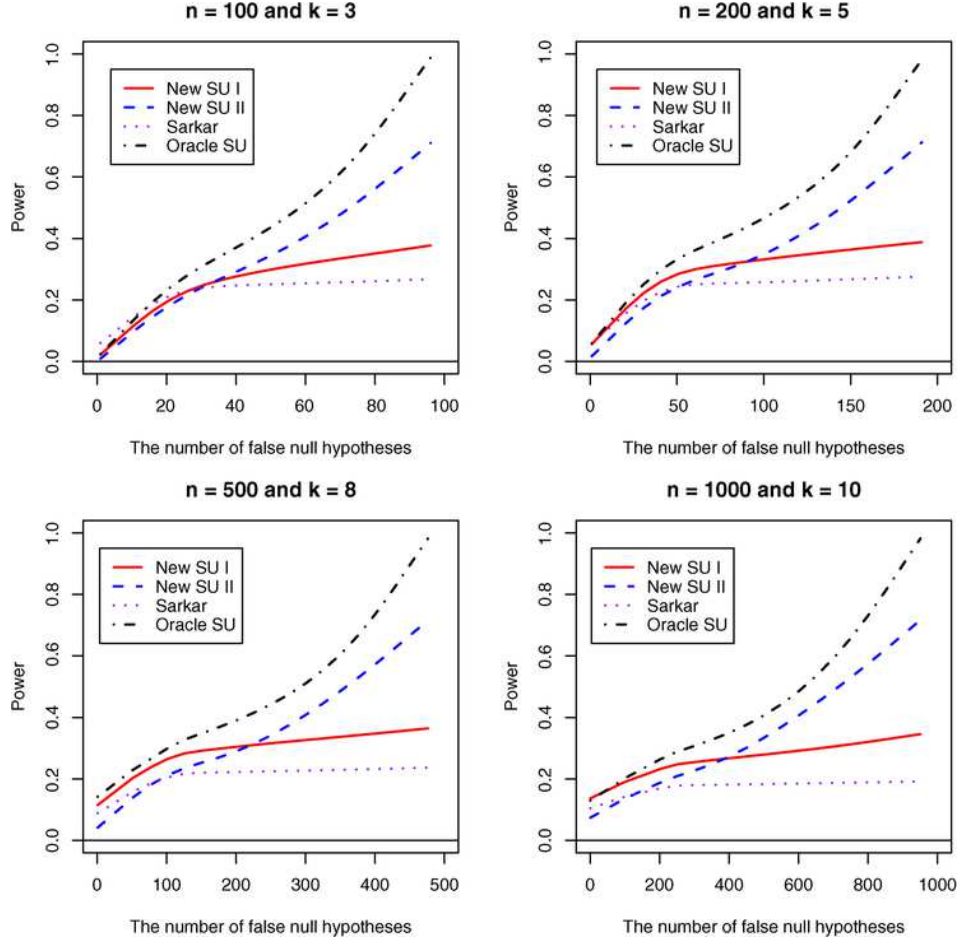

Fig. 3. *Comparison of average power of $k$-FDR stepup procedures based on the sets of critical values given by (4.2), (4.8), (5.3) and (5.4) for $\alpha = 0.05$ and $\lambda = 0.5$.*

For each gene, the base 2 logarithmic transformation of the ratio was performed before computing its two-sample $t$-test statistic. We then computed its associated raw $p$-value by using a permutation method from [25] with the permutation number $B = 1000$. Finally, we adjusted these 3170 raw $p$-values using the following five different procedures: the three $k$-FDR procedures, Procedures 1 and 2, Sarkar's $k$-FDR procedure and two $k$-FWER procedures, Sarkar's $k$-FWER and the generalized Hochberg procedures, which are now labeled in Table 1 New SU I, New SU II, SK $k$-FDR SU, SK $k$-

TABLE 1
*Numbers of differentially expressed genes for the data in [7] with $\alpha = 0.05$ and $\lambda = 0.9$*

|              | $k=1$ | $k=3$ | $k=5$ | $k=8$ | $k=10$ | $k=15$ | $k=20$ | $k=30$ |
|--------------|-------|-------|-------|-------|--------|--------|--------|--------|
| New SU I     | 74    | 75    | 81    | 103   | 124    | 157    | 173    | 229    |
| New SU II    | 129   | 129   | 129   | 135   | 137    | 162    | 176    | 229    |
| SK $k$-FDR SU  | 74    | 33    | 50    | 73    | 76     | 94     | 114    | 145    |
| SK $k$-FWER SU | 2     | 19    | 33    | 56    | 73     | 87     | 107    | 138    |
| GH $k$-FWER SU | 2     | 5     | 8     | 11    | 17     | 21     | 24     | 33     |

FWER SU and GH $k$-FWER SU, respectively. For $\alpha = 0.05$ and $\lambda = 0.9$, the numbers of significant genes found by the above five methods are presented in this table for different values of $k = 1, 3, 5, 8, 10, 15, 20$ and 30.

As expected, the $k$-FWER procedures are seen to be extremely conservative, unless $k$ is chosen large, compared to the $k$-FDR procedures. Among the $k$-FDR procedures, the two proposed ones, particularly Procedure 2, always detect much more differentially expressed genes. The SK $k$-FDR SU does not appear to be more powerful than the original BH FDR procedure unless $k$ is large (relative to $n$); whereas, those proposed here are more powerful for all $k$ (see Section 7 for further remarks on this).

**7. Concluding remarks and additional numerical investigations.** There is currently a growing interest in developing theory and methodology of multiple testing when the control of at least $k$ false rejections, for some fixed $k > 1$, rather than at least one, is of importance. A number of related procedures have been put forward in the literature, most of which are developed generalizing the traditional FWER. The generalized notion of FDR, the $k$-FDR, introduced recently in [19], on the other hand, provides a more powerful framework in this context. This is a key point, though highlighted before in [19], is re-emphasized in this paper through alternative procedures controlling the $k$-FDR, at least under independence.

A procedure controlling $k$ false discoveries should get more powerful as $k$ increases, as more and more rejections are being allowed by increasing $k$. The $k$-FDR procedures proposed here have this feature, whereas the one previously proposed in [19] does not have (as seen in Table 1). Of course, one should keep in mind that the procedure in [19] was originally developed not for independent $p$-values but for dependent $p$-values explicitly utilizing the $k$th order joint distribution of the null $p$-values. Having said that, we must nevertheless emphasize the point that even though our procedures are uniformly more powerful than the corresponding FDR procedures, one should not judge the performance of a $k$-FDR procedure against FDR procedures in the context of controlling $k$ false rejections. It should be judged, as noted in the Introduction, against compatible $k$-FWER and other $k$-FDR procedures. In fact, the difference between $k$-FDR and the corresponding FDR procedures diminishes as $n$ becomes large relative to $k$.

Choosing the value of $k$ in a $k$-FDR procedure is an important issue. It could be pre-determined. For instance, in a microarray experiment involving thousands of genes where the scientist knows that the chance of more than one hypothesis being falsely rejected is high, he/she may find it worthwhile to make further investigative studies once at least a given number, more than just one, are found differentially expressed. It could also be data-driven in that a reasonable choice of $k$ can be made only after looking at the data. For example, suppose that we are testing 100,000 hypotheses using a method controlling the FDR at 5% level. If 100 hypotheses are rejected, then one might feel comfortable adjusting this procedure to one that allows a few false rejections, say at most 9, and controls 10 or more of those at this level in an attempt to improve the power of detecting more truly false null hypotheses. On the other hand, if only 12 hypotheses are rejected, 10 is clearly not a comfortable choice. In any event, the choice of $k$ should make it more worthwhile to control the $k$-FDR than the FDR. Let us suppose that the $p$-values are independent and one likes to use our Procedure 1 to control the $k$-FDR. Notice that in this procedure the $i$th critical value $i\alpha/n$ of the BH procedure is calibrated to the $\alpha_i$ satisfying $\alpha_i G(k-1, n-1, \alpha_i) = i\alpha/n$. Thus, we have larger rejection thresholds and hence more power, and the factor $G$ essentially gives an idea about the choice of $k$ relative $n$. Let $k/n \to \gamma \in (0,1)$ as $n \to \infty$. Any $\gamma > 0$ gives more power to Procedure 1 compared to the BH procedure, but the gain in power is negligible as $\gamma \to 0$. An appropriate value of $\gamma$ can be determined subject to a desirable amount of improvement over the BH procedure. But, we will attempt to address it more formally in a different communication.

It would be interesting to see how different $k$-FDR procedures, including the oracle, proposed here under the independence assumption continue to perform in dependence cases. We did some additional simulations to investigate this. Among different possible types of dependence, we considered the equal correlation case. In particular, we generated 500 dependent normal random variables with the same variance 1 and a common correlation $\rho$, performed a multiple test using each of Procedures 1, 2 (with $\lambda = 0.5$), the oracle and the BH procedure to test each mean at 0 against 2 using a two-sided test, and repeated this over 1000 runs to simulate the $k$-FDR. The BH procedure was included for reference. Figure 4 compares the simulated $k$-FDR of these procedures with $k = 8$ and some small values of $\rho$, with Procedures 1 and 2 labeled, respectively, New SU I and New SU II. Interestingly, while both Procedure 1 and the oracle lose control of the $k$-FDR with increasing dependence, Procedure 2 seems to hold it under dependency, at least when the dependence is not too high.

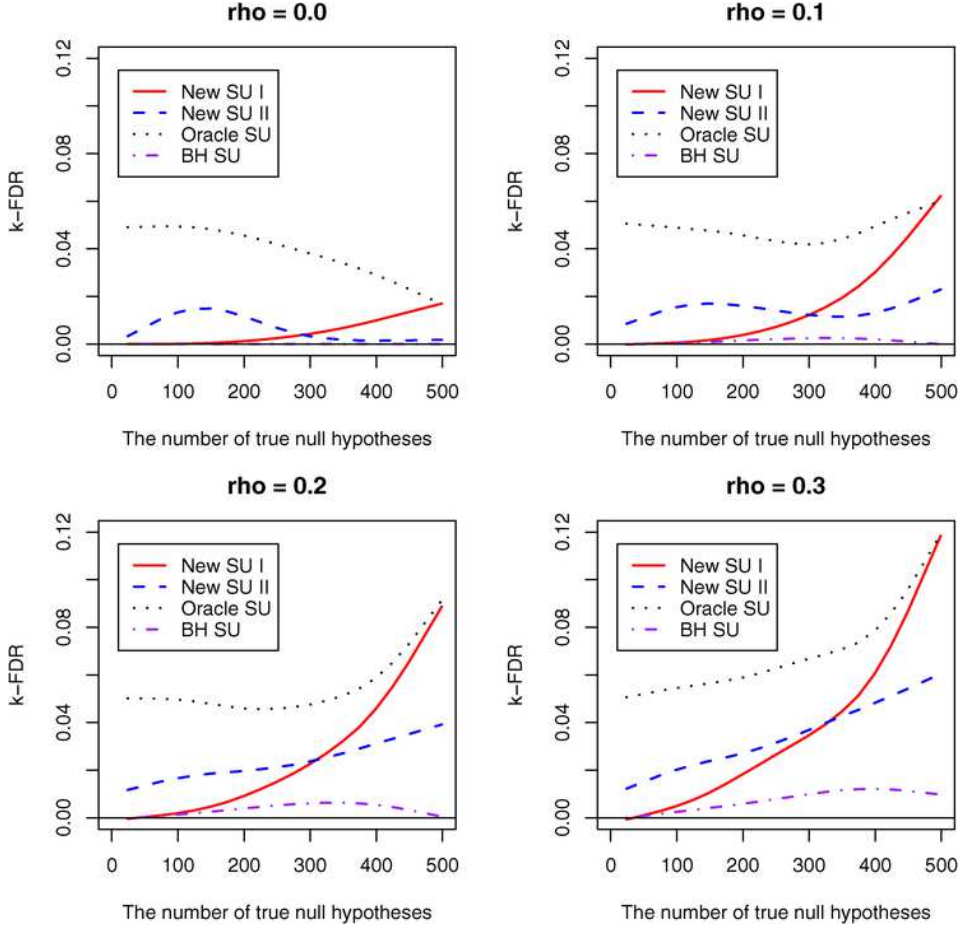

FIG. 4. *Comparison of k-FDR under dependence for $\alpha = 0.05$.*

# APPENDIX: PROOFS

**A.1. Proof of Theorem 2.1.** Define

$$V_{n-1}^{(-i)}(t) = \sum_{j(\neq i)=1}^{n} I(p_j \leq t, H_j = 0), \qquad R_{n-1}^{(-i)}(t) = \sum_{j(\neq i)=1}^{n} I(p_j \leq t).$$

Then, we note that

$$\begin{aligned}
k\text{-FDR}(t) &= E\left\{\frac{V_n(t)}{R_n(t)} I[V_n(t) \geq k]\right\} \\
&= E\left\{\sum_{i=1}^{n} \sum_{r=k}^{n} \frac{1}{r} I(p_i \leq t, H_i = 0)\right.
\end{aligned}$$

$$\times I[V_{n-1}^{(-i)}(t) + I(p_i \leq t, H_i = 0) \geq k,$$
$$R_{n-1}^{(-i)}(t) + I(p_i \leq t) = r]\bigg\}$$

(A.1)
$$= E\bigg\{\sum_{i=1}^{n}\sum_{r=k}^{n}\frac{1}{r}I(p_i \leq t, H_i = 0)$$
$$\times \Pr[V_{n-1}^{(-i)}(t) \geq k - I(p_i \leq t, H_i = 0),$$
$$R_{n-1}^{(-i)}(t) = r - I(p_i \leq t)]\bigg\}$$
$$= \pi_0 t \sum_{i=1}^{n}\sum_{r=k}^{n}\frac{1}{r}\Pr\{V_{n-1}^{(-i)}(t) \geq k-1, R_{n-1}^{(-i)}(t) = r-1\}$$
$$= n\pi_0 t \sum_{r=k}^{n}\frac{1}{r}\Pr\{V_{n-1}(t) \geq k-1, R_{n-1}(t) = r-1\}$$
$$= n\pi_0 t E\bigg\{\frac{I[V_{n-1}(t) \geq k-1]}{R_{n-1}(t) + 1}\bigg\}.$$

The probability in the third line in (A.1) is obtained by taking the conditional expectation given $I(p_i \leq t)$ and $I(H_i = 0)$ of the inner indicator function in the previous line. The fourth line follows from the fact that the expectation of the product of $I(p_i \leq t, H_i = 0)$ and a function of $I(p_i \leq t)$ and $I(H_i = 0)$ is $\Pr\{p_i \leq t, H_i = 0\}$ times the value of that function when both $I(p_i \leq t)$ and $I(H_i = 0)$ are 1.

**A.2. Proof of Theorem 4.1.** Let us first prove the following two results that will be useful in proving the theorem.

RESULT A.1. *The function $G(k, n, u)$ defined in (2.10) is nondecreasing in $n$ and $u$, for any fixed $1 \leq k \leq n$.*

PROOF. Note that $G(k, n, u) = \Pr(U_{k:n} \leq u)$ with $U_{k:n}$ being the $k$th order statistic based on $n$ i.i.d. Uniform$(0,1)$ random variables, which is clearly increasing in $u$ for fixed $k$ and $n$. Since the value of $U_{k:n}$ decreases as $n$ increases, $G$ is also increasing in $n$ for fixed $k$ and $u$. □

RESULT A.2. *Let $R_n$ be the total number of rejections in a stepup procedure based on p-values $p_1, \ldots, p_n$ and critical values $\alpha_1 \leq \cdots \leq \alpha_n$. Let $\hat{p}_{1:n_0} \leq \cdots \leq \hat{p}_{n_0:n_0}$ be the ordered p-values corresponding to the $n_0$ true null* *hypotheses. Then, for any fixed* $1 \leq k \leq n_0$,

$$(A.2) \qquad \sum_{r=k}^{n} \Pr(R_n = r | \hat{p}_{k:n_0} \leq \alpha_r) \leq 1.$$

PROOF.

$$\sum_{r=k}^{n} \Pr(R_n = r | \hat{p}_{k:n_0} \leq \alpha_r)$$

$$= \sum_{r=k}^{n} \Pr(R_n \geq r | \hat{p}_{k:n_0} \leq \alpha_r) - \sum_{r=k}^{n-1} \Pr(R_n \geq r+1 | \hat{p}_{k:n_0} \leq \alpha_r)$$

$$= \Pr(R_n \geq k | \hat{p}_{k:n_0} \leq \alpha_k)$$

$$+ \sum_{r=k}^{n-1} [\Pr(R_n \geq r+1 | \hat{p}_{k:n_0} \leq \alpha_{r+1}) - \Pr(R_n \geq r+1 | \hat{p}_{k:n_0} \leq \alpha_r)].$$

The result then follows from the fact that $\Pr(R_n \geq k | \hat{p}_{k:n_0} \leq \alpha_k) = 1$, because, given the occurrence of at least $k$ false rejections, the probability that at least $k$ hypotheses are rejected is 1, and that

$$(A.3) \qquad \Pr(R_n \geq r+1 | \hat{p}_{k:n_0} \leq \alpha_{r+1}) \leq \Pr(R_n \geq r+1 | \hat{p}_{k:n_0} \leq \alpha_r),$$

for all $r = k, \ldots, n-1$, which can be proved as follows.

Since $R_n$ decreases as each of $\hat{p}_{1:n_0}, \ldots, \hat{p}_{n_0:n_0}$ and the nonnull $p$-values increases, the conditional probability

$$(A.4) \qquad g(\hat{p}_{1:n_0}, \ldots, \hat{p}_{n_0:n_0}) = \Pr(R_n \geq r+1 | \hat{p}_{1:n_0}, \ldots, \hat{p}_{n_0:n_0})$$

is nonincreasing in $\hat{p}_{1:n_0}, \ldots, \hat{p}_{n_0:n_0}$. Now, the order statistics, say $X_{1:m} \leq \cdots \leq X_{m:m}$, of any set of $m$ i.i.d. (continuous) random variables are stochastically increasing in each of its components, that is, $E\{\phi(X_{1:m}, \ldots, X_{m:m}) | X_{k:m}\}$ is nondecreasing (or nonincreasing) in $X_{k:m}$, for any fixed $1 \leq k \leq m$ and nondecreasing (or nonincreasing) function $\phi$. See, for example, Block, Savits and Shaked [4], who defined this condition as the positive dependent through stochastic ordering (PDS) condition. It also follows from the positive regression dependence condition satisfied by the joint distribution of order statistics; see, for example, Karlin and Rinott [11]. Thus, the conditional probability

$$(A.5) \qquad \begin{aligned} h(\hat{p}_{k:n_0}) &= \Pr(R_n \geq r+1 | \hat{p}_{k:n_0}) \\ &= E\{g(\hat{p}_{1:n_0}, \ldots, \hat{p}_{n_0:n_0}) | \hat{p}_{k:n_0}\}, \end{aligned}$$

is nonincreasing in $\hat{p}_{k:n_0}$, and hence

$$(A.6) \qquad \Pr(R_n \geq r+1 | \hat{p}_{k:n_0} \leq t) = \frac{E\{h(\hat{p}_{k:n_0}) I(\hat{p}_{k:n_0} \leq t)\}}{\Pr(\hat{p}_{k:n_0} \leq t)},$$

is nonincreasing in $t$. $\square$

We are now ready to prove the theorem. First, note that given $n_0$ true null hypotheses with the corresponding $p$-values $\hat{p}_1, \ldots, \hat{p}_{n_0}$, the $k$-FDR of a stepup procedure with critical values $\alpha_1 \leq \cdots \leq \alpha_n$ under the conditions assumed in the theorem is

$$\text{(A.7)} \qquad k\text{-FDR} = n_0 \sum_{r=k}^{n} \frac{\alpha_r}{r} \Pr(V_{n-1} \geq k-1, R_{n-1} = r-1),$$

where $R_{n-1}$ and $V_{n-1}$ are the number of rejections and the number of false rejections, respectively, in the stepup procedure based on the $n-1$ $p$-values $\{p_1, \ldots, p_n\} \setminus \{\hat{p}_{n_0}\}$ and the critical values $\alpha_i$, $i = 2, \ldots, n$. Let $\hat{p}_{1:n_0-1} \leq \cdots \leq \hat{p}_{n_0-1:n_0-1}$ be the ordered $n_0 - 1$ null $p$-values. Then, given $\{R_{n-1} = r-1\}$, $V_{n-1} \geq k-1$ if and only if $\hat{p}_{k-1:n_0-1} \leq \alpha_r$. Thus, the $k$-FDR in (A.7) is equal to

$$n_0 \sum_{r=k}^{n} \frac{\alpha_r}{r} \Pr(\hat{p}_{k-1:n_0-1} \leq \alpha_r) \Pr(R_{n-1} = r-1 | \hat{p}_{k-1:n_0-1} \leq \alpha_r)$$

$$\text{(A.8)} \qquad = n_0 \sum_{r=k}^{n} \frac{\alpha_r}{r} G(k-1, n_0-1, \alpha_r) \Pr(R_{n-1} = r-1 | \hat{p}_{k-1:n_0-1} \leq \alpha_r)$$

$$\leq n_0 \sum_{r=k}^{n} \frac{\alpha_r}{r} G(k-1, n-1, \alpha_r) \Pr(R_{n-1} = r-1 | \hat{p}_{k-1:n_0-1} \leq \alpha_r),$$

using Result A.1. For the stepup procedure (4.2) with its critical values satisfying $\alpha_r G(k-1, n-1, \alpha_r) = r\alpha/n$, $r = 1, \ldots, n$, that are increasing because of Result A.1, we have

$$\text{(A.9)} \qquad k\text{-FDR} \leq \frac{n_0}{n} \alpha \sum_{r=k}^{n} \Pr(R_{n-1} = r-1 | \hat{p}_{k-1:n_0-1} \leq \alpha_r).$$

The theorem then follows form Result A.2.

**A.3. Proof of Theorem 4.2.** Our proof relies on arguments used in proving Theorem 4.1 and in [18]. Define

$$\text{(A.10)} \quad \alpha_{ij} = \lambda \min\left\{ \tilde{G}_{k,n}^{-1}\left( \frac{i\alpha(1-\lambda)}{\lambda(n-j+1)} \right), 1 \right\}, \qquad 1 \leq i \leq j; \; j = 1, \ldots, n.$$

Consider $E(k\text{-FDP}|p_{j:n} \leq \lambda < p_{j+1:n})$, the $k$-FDR conditional on $p_{j:n} \leq \lambda < p_{j+1:n}$. This conditional $k$-FDR is the $k$-FDR of the stepup procedure based on $j$ independent $p$-values, each truncated above at $\lambda$ and the critical values $\alpha_{ij}$, $i = 1, \ldots, j$. Noting that these truncated $p$-values corresponding to the true null hypotheses are i.i.d. Uniform$(0, \lambda)$ and that

$$\frac{\alpha_{rj}}{\lambda r} G\left(k-1, n-1, \frac{\alpha_{rj}}{\lambda}\right) = \frac{1}{r} \tilde{G}_k\left(\frac{\alpha_{rj}}{\lambda}\right) \leq \frac{\alpha(1-\lambda)}{\lambda(n-j+1)},$$

we see by arguing as in our proof of Theorem 4.1 [see (A.8)] that this conditional $k$-FDR is less than or equal to $(1-\lambda)\alpha/[\lambda(n-j+1)]E\{V(\lambda)|p_{j:n} \leq \lambda < p_{j+1:n}\}$, where $V(\lambda) = \sum_{i=1}^{n_0} I(\hat{p}_i \leq \lambda)$. Thus, we have

$$k\text{-FDR} = \sum_{j=k}^{n} E(k\text{-FDP}|p_{j:n} \leq \lambda < p_{j+1:n}) \Pr(p_{j:n} \leq \lambda < p_{j+1:n})$$

$$\leq \alpha E\left\{\sum_{j=k}^{n} \frac{1-\lambda}{\lambda(n-j+1)} V(\lambda) I(p_{j:n} \leq \lambda < p_{j+1:n})\right\}$$

$$= \alpha E\left\{\sum_{i=1}^{n_0} \sum_{j=k}^{n} \frac{1-\lambda}{\lambda(n-j+1)} I(\hat{p}_i \leq \lambda, p_{j:n} \leq \lambda < p_{j+1:n})\right\}$$

(A.11)
$$\leq n_0 \alpha \sum_{j=k}^{n} \frac{1-\lambda}{n-j+1} \Pr(p_{j-1:n-1} \leq \lambda < p_{j:n-1})$$

$$= \alpha \sum_{i=1}^{n_0} \sum_{j=k}^{n} \frac{1}{n-j+1} \Pr(\hat{p}_i > \lambda) \Pr(p_{j-1:n-1} \leq \lambda < p_{j:n-1})$$

$$\leq \alpha \sum_{i=1}^{n} \sum_{j=k}^{n} \frac{1}{n-j+1} \Pr(p_i > \lambda) \Pr(p_{j-1:n-1} \leq \lambda < p_{j:n-1})$$

$$\leq \alpha \sum_{j=k}^{n} \Pr(p_{j-1:n} \leq \lambda < p_{j:n}) = \alpha \Pr(p_{k-1:n} \leq \lambda < p_{n:n}) \leq \alpha,$$

proving the theorem.

**A.4. Proof of (2.7).** Using the result that two functions, one increasing and the other decreasing, of a random variable are negatively correlated, we first see that

(A.12)
$$E\left\{\frac{I(V_{n-1} \geq k-1)}{R_{n-1}(t)+1}\right\} = E\left\{\frac{\Pr(V_{n-1} \geq k-1|R_{n-1})}{R_{n-1}(t)+1}\right\}$$

$$\leq E\left\{\frac{1}{R_{n-1}(t)+1}\right\} \Pr\{V_{n-1} \geq k-1\}.$$

The inequality then follows by dividing both sides of (A.12) by the probability $\Pr\{V_{n-1} \geq k-1\}$.

**Acknowledgments.** We are grateful to two referees whose insightful comments have considerably improved the paper.

DEPARTMENT OF STATISTICS
FOX SCHOOL OF BUSINESS AND MANAGEMENT
TEMPLE UNIVERSITY
PHILADELPHIA, PENNSYLVANIA 19122
USA
E-MAIL: sanat@temple.edu

NATIONAL INSTITUTE OF ENVIRONMENTAL
  HEALTH SCIENCES
BIOSTATISTICS BRANCH
PO BOX 12233 MAIL DROP A3-03
111 TW ALEXANDER DRIVE
RESEARCH TRIANGLE PARK,
  NORTH CAROLINA 27709-2233
USA
E-MAIL: guow@niehs.nih.gov